\documentclass[12]{amsart}
\usepackage{latexsym}
\usepackage{amssymb}
\usepackage{amsfonts, amsmath}
\usepackage[T1]{fontenc}
\usepackage[utf8]{inputenc}
\usepackage{newunicodechar}
\newunicodechar{−}{-}
\usepackage{amsmath, amssymb, amsthm, latexsym}

\newtheorem{lm}{Lemma}
\newtheorem{Thm}{Theorem}

\newtheorem*{corr}{Corollary}

\begin{document}
\title{The Local Steiness problem with singularities.}
\author{}
\noindent
\maketitle
\begin{center}
\em{Youssef Alaoui}\\
\em{y.alaoui@iav.ac.ma}\\
\end{center}
\noindent
\em{Department of Mathematics,
Hassan II Institute of Agronomy}\\ 
\em{and Veterinary Sciences,
Madinat Al Irfane, BP 6202, Rabat, 10101, Morocco,}\\
\date{}
\linespread{1.3}
\date{}
\begin{abstract}
In this article, we prove that if $\Pi: X\rightarrow
\Omega$ is an unbranched Riemann domain with $\Omega$
Stein of dimension $n$ and $\Pi$ a locally $q$-complete morphism, then $X$
is cohomologically $q$-complete if $n\geq 3$ and $1\leq q\leq n-2$ or if $\Omega$ has dimension $2$
and $1\leq q\leq 2$.
This generalizes a well-known result which is obtained in ~\cite{ref3} for $q=1$
when $X$ and $\Omega$ have isolated singularities and, gives in particular a positive
answer to the local Steiness problem, namely if $X$ is a Stein space and $\Omega$
a locally Stein open subset of $X$, then $\Omega$ is Stein.
\end{abstract}
\maketitle
\section{Introduction}
A Riemann domain over a complex space $Y$ is a pair $(X, \Pi)$, where
$\Pi: X\rightarrow Y$ is a holomorphic map which is non-degenerate at every point
of $X$, i.e. $\Pi^{-1}(\Pi(x))$ is a discreate set at each point $x\in X$.
The pair $(X, \Pi)$ is called unbranched or unramified
if $\Pi: X\rightarrow Y$ is locally biholomorphic.\\
\hspace*{.1in}A Riemann domain $(X,\Pi)$ over $Y$ is locally $q$-complete
if  there exists for every $x\in Y$ an open neighborhood $U$ in $Y$
such that $\Pi^{-1}(U)$ is $q$-complete.\\
\hspace*{.1in}Let $X$ and $Y$ be complex spaces and $\Pi: X\rightarrow Y$ an unbranched Riemann domain such that $Y$ is Stein and $\Pi$ a locally $q$-complete morphism.\\
Does it follow that $X$ is $q$-complete ?\\
\hspace*{.1in}In ~\cite{ref3}, it was proved that this problem has a positive answer if $q=1$ and $X$ and $Y$ have isolated singularities.\\
\hspace*{.1in}It has been proved however in ~\cite{ref5} that there exist locally
$1$-complete Riemann domains over $\mathbb{C}^{n}$
which are not holomorphically-convex.\\
\hspace*{.1in}It is known from ~\cite{ref7} that if $\Pi : X\rightarrow \Omega$ is an unbranched Riemann domain between two complex spaces with isolated singularities, $\Omega$ $q$-complete and $\Pi$ is locally $1$-complete, then $X$ is $q$-complete.\\
\hspace*{.1in}In this article, we prove that if
$\Pi : X\rightarrow \Omega$ is a locally $q$-complete unbranched Riemann domain
over a Stein space $\Omega$ of dimension $n$, then for any
coherent analytic sheaf ${\mathcal{F}}$ on $X$,
the cohomology group $H^{q}(X,{\mathcal{F}})=0$, if $n\geq 3$ and $1\leq q\leq n-2$ or if $dim \Omega$
has dimension $2$ and $1\leq q\leq 2$.\\
\hspace*{.1in}In particular, we obtain the interesting result
\begin{corr}
{If $X$ is a Stein space of dimension $n\geq 3$ and $\Omega\subset X$ a locally $q$-complete
open subset of $X$, then $\Omega$ is cohomologically $q$-complete if
$1\leq q\leq n-2$ or if $X$ has dimension $2$ and $1\leq q\leq 2$.}
\end{corr}
It is well known from ~\cite{ref11} that if $Y$ is Stein and, if  $\Pi: X\rightarrow Y$ is a locally $q$-complete morphism, then $X$ is cohomologically $(q+1)$-complete.
But in general $H^{q}(X, {\mathcal{O}}_{X})$ does not vanish,
even when $\Pi: X\rightarrow Y$ is locally Stein ~\cite{ref9}
(See also ~\cite{ref4}).
\section{Preliminaries}
We start by recalling some definitions and results concerning $q$-complete spaces.\\
\hspace*{.1in}Let $\Omega$ be an open set in $\mathbb{C}^{n}$
with complex coordinates $z_{1},\cdots, z_{n}$. Then it is known that
a function $\phi\in C^{\infty}(\Omega)$ is $q$-convex if for every
point $z\in \Omega$, the Levi form
$$L_{z}(\phi,\xi)=\displaystyle\sum_{i,j}\frac{\partial^{2}\phi(z)}{\partial{z_{i}}\partial{\overline{z}_{j}}}\xi_{i}\overline{\xi}_{j}, \ \ \xi\in \mathbb{C}^{n}$$
has at most $q-1$ negative or zero eingenvalues.\\
\hspace*{.1in}A smooth real valued function $\phi$ on a complex space $X$
is called $q$-convex if every point $x\in X$ has a local chart
$U\rightarrow D\subset \mathbb{C}^{n}$ such that $\phi|_{U}$ has an extension $\hat{\phi}\in C^{\infty}(D, \mathbb{R})$ which is $q$-convex on $D$.\\
\hspace*{.1in}We say that $X$ is $q$-complete if there exists
a $q$-convex function $\phi\in C^{\infty}(X,\mathbb{R})$
which is exhaustive on $X$ i.e. $\{x\in X: \phi(x)<c\}$
is relatively compact for any $c\in\mathbb{R}.$\\
\hspace*{.1in}The space $X$ is said to be cohomologically $q$-complete
if for every coherent analytic sheaf ${\mathcal{F}}$ on $X$ the cohomology groups $H^{p}(X, {\mathcal{F}})$ vanish for all $p\geq q$.\\
\hspace*{.1in}An open subset $D$ of $\Omega$ is called $q$-Runge
if for every compact set $K\subset D$, there is a $q$-convex
exhaustion function $\phi\in C^{\infty}(\Omega)$ such that
$$K\subset\{x\in\Omega: \phi(x)<0\}\subset\subset D$$
This generalizes the classical notion of Runge pairs of Stein spaces.\\
It is shown in $[3]$ that if $D$ is $q$-Runge in $\Omega$, then for
every ${\mathcal{F}}\in coh(\Omega)$ the cohomology groups $H^{p}(D,
\mathcal{F})$ vanish for $p\geq q$ and, the restriction map
$$H^{p}(\Omega, {\mathcal{F}})\longrightarrow H^{p}(D, {\mathcal{F}})$$
has dense image for all $p\geq q-1$.\\
\begin{lm}{Let $X$ and $Y$ be complex manifolds of dimension $n$ and $\Pi:
X\rightarrow Y$ an unbranched Riemann domain. Assume that there
exists a smooth $q$-convex function $\phi$ on $Y$. Let $\xi_{0}\in X$
and $X'_{c}=\{x\in X: \phi o\Pi(x)>c\},$ where $c=\phi o\Pi(\xi_{0})$.
Then for any coherent analytic sheaf ${\mathcal{F}}$
on $X$ the restriction map
$$H^{p}(X,{\mathcal{F}})\rightarrow H^{p}(X'_{c}, {\mathcal{F}})$$
is bijective if $p\leq n-q-1$,\\
\hspace*{.1in}injective if $p=n-q.$}
\end{lm}
Let $D$ be a domain in $\mathbb{C}^{n}$, $\xi\in D$,
and let $\phi\in C^{\infty}(D )$ be a q-convex
function. Then in order to prove lemma $1$ we shall need the following result due to Andreotti and Grauert ~\cite{ref2}.
\begin{Thm}{For any coherent analytic sheaf ${\mathcal{F}}$ on $D$
there exists a fundamental system of Stein
neighborhoods $U\subset D$ of $\xi$ such that if $Y=\{z\in D:
\phi(z)>\phi(\xi)\}$, then $H^{p}(Y\cap U, {\mathcal{F}})=0$ for
$0<p<n-q$ and $H^{0}(U,{\mathcal{F}})\rightarrow H^{0}(U\cap Y, {\mathcal{F}})$ is an isomorphism.}
\end{Thm}
\noindent
\begin{proof}
Let $V\subset\subset X$ be a hyperconvex open neighborhood of $\xi_{0},$
biholomorphic by $\Pi$ to the open subset $W=\Pi(V)\subset Y$. We
may take $V$ so that $W$ is biholomorphic to a domain $D$ in $\mathbb{C}^{n}$.\\
\hspace*{.1in} Let $\psi: V\rightarrow ]-\infty, 0[$ be a continuous strictly plurisubharmonic function. Then it is clear that $\psi_{k}=\frac{1}{k}\psi+\phi o\Pi,$ $k\geq
1,$ is an increasing sequence of $q$-convex functions on $V.$ If we
put $V_{k}=\{x\in V: \psi_{k}(x)>c\},$ then
$\displaystyle\bigcup_{k\geq 1}V_{k}=V\cap X'_{c}.$ Moreover, there exists, by theorem $1,$ a fundamental system of connected Stein neighborhoods $U\subset V$ of $\xi_{0}$ such that $H^{r}(U\cap V_{k}, {\mathcal{F}})=0$ for $1\leq r<n-q$ and
$H^{0}(U, {\mathcal{F}})\rightarrow H^{0}(U\cap V_{k},
{\mathcal{F}})$ is an isomorphism, or equivalently (See ~\cite{ref7} or
~\cite{ref1}), $\underline{H^{r}_{S_{k}}}({\mathcal{F}})=0$ for $r\leq
n-q,$ where
$\underline{H^{r}_{S_{k}}}({\mathcal{F}})$ is the cohomology sheaf
with support in $S_{k}=\{x\in V: \psi_{k}(x)\leq c\}$ and
coefficients in ${\mathcal{F}}.$ Furthermore, there exists a
spectral sequence
$$H^{p}_{S_{k}}(V, {\mathcal{F}})\Longleftarrow
E_{2}^{p,q}=H^{p}(V, \underline{H^{q}_{S_{k}}}({\mathcal{F}}))$$
Since $\underline{H^{p}_{S_{k}}}({\mathcal{F}})=0$ for $p\leq n-q,$ then for any $p\leq n-q,$
the cohomology groups $H^{p}_{S_{k}}(V, {\mathcal{F}})$ vanish and, the exact sequence of local
cohomology
\begin{center}
$\cdots\rightarrow H^{p}_{S_{k}}(V, {\mathcal{F}}) \rightarrow
H^{p}(V, {\mathcal{F}}) \rightarrow H^{p}(V_{k}, {\mathcal{F}})
\rightarrow H^{p+1}_{S_{k}}(V, {\mathcal{F}})\rightarrow\cdots$
\end{center}
implies that $H^{p}(V_{k}, {\mathcal{F}})\cong H^{p}(V,
{\mathcal{F}})$ for all $p\leq n-q-1.$ Hence\\
$H^{p}(V_{k}, {\mathcal{F}})=0$ for $1\leq p\leq
n-q-1$ and $H^{o}(V_{k}, {\mathcal{F}})\cong
H^{o}(V, {\mathcal{F}})$ for every integer $k$. Since $V\cap X'_{c}$
is an increasing union of $V_{k}$, $k\in \mathbb{N}$, then, by (~\cite{ref2},
lemma, p.$250$), we deduce that $H^{p}(V\cap
X'_{c},{\mathcal{F}})=0$ for $1\leq p\leq n-q-1$ and $H^{o}(V,
{\mathcal{F}})\rightarrow H^{o}(V\cap X'_{c}, {\mathcal{F}})$ is an
isomorphism. Since each point of $X$ has a
fundamental system of hyperconvex neighborhoods, then, if
$S=\{x\in X: \phi o\Pi(x)\leq c\},$ the cohomology sheaf
$\underline{H^{p}_{S}}({\mathcal{F}})$ vanishes for all $p\leq
n-q.$ Therefore the spectral sequence
$$H^{p}_{S}(X, {\mathcal{F}})\Longleftarrow
E_{2}^{p,q}=H^{p}(X, \underline{H^{q}_{S}}({\mathcal{F}}))$$ shows
that $H^{p}_{S}(X, {\mathcal{F}})=0$ for any $p\leq
n-q,$ and from the exact sequence
$$\cdots\rightarrow H^{p}_{S}(X, {\mathcal{F}})\rightarrow H^{p}(X, {\mathcal{F}})\rightarrow H^{p}(X'_{c}, {\mathcal{F}})\rightarrow H^{p+1}_{S}(X, {\mathcal{F}})\rightarrow\cdots$$
we see that the restriction map $H^{r}(X, {\mathcal{F}})\rightarrow H^{r}(X'_{c},{\mathcal{F}})$ is bijective
for any $c\in\mathbb{R}$ if $r\leq n-q-1$ and,
injective if $r=n-q$.
\end{proof}
\begin{lm}Let $X$ and $Y$ be two $n$-dimensional complex manifolds such that
$Y$ is Stein and $\Pi: X\rightarrow Y$ is an unbranched Riemann domain
and locally $q$-complete morphism with $1\leq q\leq n-2$. Then $X$ is cohomologically $q$-complete.
\end{lm}
\begin{proof}
\hspace*{.1in}In fact, we consider a covering ${\mathcal{V}} = (V_{i})_{i\in\mathbb{N}}$ of $Y$ by open sets $V_{i}\subset Y$ such that
$\Pi^{-1}(V_{i})$ is $q$-complete for all $i\in\mathbb{N}$. By the Stein covering lemma of Sthel\'e ~\cite{ref10}, there exists a locally finite covering ${\mathcal{U}} = (U_{i})_{i\in\mathbb{N}}$ of $Y$ by Stein open
subsets $U_{i}\subset\subset Y$ such that ${\mathcal{U}}$ is a refinement of $ {\mathcal{V}}$,
$\displaystyle\bigcup_{i\leq j}U_{i}$ is Stein for all is Stein for all j.
Moreover, there exists for all $j\in\mathbb{N}$ a continuous strictly plurisubharmonic
function $\phi_{j+1}$ on $\displaystyle\bigcup_{i\leq j+1}U_{i}$ such that
$$\displaystyle\bigcup_{i\leq j}U_{i}\cap U_{j+1} =\{ x \in U_{j+1} : \phi_{j+1}(x) < 0\}$$
Note also that $\Pi^{-1}(U_{i})$ is $q$-complete for all $i\in\mathbb{N}$ and, if
$X_{j} = \Pi^{-1}(\displaystyle\bigcup_{i\leq j}U_{i})$
and $X'_{j+1}=\Pi^{-1}(U_{j+1})$, then $X_{j}\cap X'_{j+1}=\{x\in X'_{j+1}: \phi_{j+1}o \Pi(x)<0\}$
is clearly $q$-Runge in $X'_{j+1}$.\\
\hspace*{.1in}Let now $\mathcal{F}$ be a coherent analytic  sheaf on $X$.
We shall first prove by induction on $j$ that the cohomology groups $H^{q}(X_{j},{\mathcal{F}})$
vanish. For $j = 0$, this is clear, since $\Pi^{-1}(U_{o})$ is $q$-complete. Suppose that $j\geq 1$,
$H^{q}(X_{j},{\mathcal{F}})=0$ and put $Y_{j}=\{x\in X_{j} : \phi_{j+1}o\Pi(x) > 0\}$ and
$Y'_{j+1}=\{x\in X'_{j+1} : \phi_{j+1}o\Pi(x) > 0\}$.\\
\noindent
\hspace*{.1in}Since $n\geq 3$ and $1\leq q\leq n-2$, then, by lemma $1$, $H^{p}(Y_{j},{\mathcal{F}})\cong H^{p}(X_{j},{\mathcal{F}})$
and $H^{p}(Y'_{j+1},{\mathcal{F}})\cong H^{p}(X'_{j+1},{\mathcal{F}})$
for $p\leq q.$ Since\\ $Y''_{j+1} =\{x\in  X_{j+1}: \phi_{j+1}o\Pi(x)>0\}=Y_{j}\cup Y'_{j+1}$ and $Y_{j}\cap Y'_{j+1}=\emptyset$, then we have
$$H^{p}(X_{j+1},{\mathcal{F}})\cong H^{p}(Y''_{j+1},{\mathcal{F}})\cong H^{p}(Y_{j},{\mathcal{F}})\oplus H^{p}(Y'_{j+1} ,{\mathcal{F}}) \ \ for \ \ all \ \ p\leq q$$
This proves in particular that $H^{q}(X_{j},{\mathcal{F}})=0$ for all $j\in\mathbb{N}.$\\
\hspace*{.1in}Moreover, since $X$ is an increasing union of $(X_{j})_{j\geq 0}$ and $H^{q-1}(X_{j+1},{\mathcal{F}})\cong
H^{q-1}(X_{j},{\mathcal{F}})\oplus H^{q-1}(X'_{j+1},{\mathcal{F}})$,
then, by (~\cite{ref2}, lemma, p. $250$),
the restriction map
$H^{q}(X,{\mathcal{F}})\rightarrow H^{q}(X_{0},{\mathcal{F}})$
is an isomorphism, which implies that $H^{q}(X,{\mathcal{F}})=0.$\\
\end{proof}
\begin{lm}Let $\Pi: X\rightarrow \Omega$ be an unbranched Riemann domain
and locally $q$-complete morphism over a normal Stein space $\Omega$ of dimension $n$. Then $X$ is cohomologically $q$-complete if $n=3$ and $1\leq q\leq n-2$ or if $n=2$ and $1\leq q\leq 2$.
\end{lm}
\begin{proof}
If $\Omega$ has dimension $2$ and $q=1$, it follows from \cite{ref3} that $X$ is Stein and,
when $q=2$, then $X$ is $2$-complete, according to a theorem of  Ohsawa ~\cite{ref8}, since obviously every compact analytic subset of $X$ is finite.\\
\hspace*{.1in}We may assume that $X$ is connected, $n=dim(\Omega)\geq 3$ and
that the theorem has already been proved in dimension $\leq n-1$.\\
Since $\Pi: X \rightarrow \Omega$ is locally $q$-complete, it follows from ~\cite{ref11} that for every coherent analytic sheaf
${\mathcal{F}}$ on $X$ the cohomology group $H^{p}(X,{\mathcal{F}})$ vanishes for all $p\geq q + 1.$
It is therefore enough to prove that $H^{q}(X,{\mathcal{F}}) = 0.$\\
\hspace*{.1in}let $f$ be a holomorphic function on $\Omega$ such that
$Sing(\Omega)\subset Z=\{x\in\Omega : f(x)=0\}$ and put $Z'=\Pi^{-1}(Z)$.
We first prove that $Z'=\Pi^{-1}(Z)$ is cohomologically $q$-complete.
In fact, let $Z_{1}\stackrel{\psi}\rightarrow Z$ be a normalization of $Z$.
If $\tilde{Z}'=\{(z,w)\in Z'\times Z_{1} : \Pi(z)=\psi(w)\},$ then the projection
$\Pi_{2} : \tilde{Z}'\rightarrow Z_{1}$ is clearly an unbrached Riemann domain
over the Stein normal space $Z_{1}$ and locally $q$-complete morphism. Therefore, by the induction hypothesis, $\tilde{Z}'$ is cohomologically $q$-complete. On the other hand, it is easy to verify that the projection $\Pi_{1} : \tilde{Z}'\rightarrow Z'$ is a finite holomorphic surjection, which implies that $Z'$ is also cohomologically $q$-complete. (See e.g. ~\cite{ref12}).\\
\hspace*{.1in}It is also clear that if $X'=X\setminus Z'$, then $H^{p}(X',{\mathcal{F}})=0$
for all $p\geq q$.\\
In fact, since the restriction map $\Pi|_{X'} : X'\rightarrow \Omega\setminus Z $
is obviously an unbranched Riemann domain and locally $q$-complete morphism over the Stein manifold $\Omega\setminus Z$,
then, according to lemma $2$, $X'$ is cohomologically $q$-complete.\\
\hspace*{.1in}Let now $\xi : \tilde{X}\rightarrow X$ be a resolution of singularities,
i.e. $\tilde{X}$ is a complex manifold and $\xi$ is a proper modification such that the induced
map
$$\tilde{X}\backslash\xi^{-1}(Sing(X))\rightarrow X\setminus Sing(X)$$
is biholomorphic. Then there exists a canonical sheaf homomorphism
${\mathcal{F}}\stackrel{\psi}\rightarrow \xi_{*}\xi^{*}({\mathcal{F}}).$
If we set ${\mathcal{F}}_{1}=ker \ \psi$ and ${\mathcal{F}}_{2}=Im \ \psi$,
then $Supp ({\mathcal{F}}_{1})\subset Z'$ and there is an exact sequence
$$0\rightarrow {\mathcal{F}}_{1}\rightarrow {\mathcal{F}}\stackrel{\psi}\rightarrow {\mathcal{F}}_{2}\rightarrow 0$$
Let ${\mathcal{I}}(Z')\subset {\mathcal{O}}_{X}$ be the subsheaf of germs of holomorphic functions which vanish on $Z'$. Let ${\mathcal{O}}_{Z'}={\mathcal{O}}_{X}/_{{\mathcal{I}}(Z')}$
and $({\mathcal{F}}_{2})_{Z'}={\mathcal{F}}_{2}\otimes_{{\mathcal{O}}_{X}}{\mathcal{O}}_{Z'}$.
If $e$ is the image in ${\mathcal{O}}_{Z'}$ of the section $1$ on ${\mathcal{O}}_{X}$,
then any element of $(({\mathcal{F}}_{2})_{Z'})_{x}$ can be written in the form
$\xi\otimes e_{x},$
where $\xi\in ({\mathcal{F}}_{2})_{x}$. Then the homomorphism
$\eta : {\mathcal{F}}_{2}\rightarrow ({\mathcal{F}}_{2})_{Z'}$ defined by
$\eta(\alpha)=\alpha\otimes e$ is surjective and we have the exact sequence
$$0\rightarrow Ker (\psi)\rightarrow Ker (\eta o\psi)\rightarrow \frac{Ker (\eta o\psi)}{Ker (\psi)}\rightarrow 0$$
Since clearly $Supp \ Ker \ \psi\subset Sing(X)\subset Z'$ and
$Supp \ \frac{Ker (\eta o\psi)}{Ker (\psi)}\subset X\setminus Z'$, it follows from
the long exact sequence of cohomology
$$\cdots\rightarrow H^{p}(Z', Ker \ \psi)\rightarrow H^{p}(X, Ker \ \eta o\psi )
\rightarrow H^{p}(X\setminus Z', \frac{Ker \ \eta o\psi}{Ker \ \psi})\rightarrow\cdots$$
that $H^{p}(X, Ker \ \eta o\psi )=0$ for all $p\geq q$.
Furthermore, since $H^{p}(X,({\mathcal{F}}_{2})_{Z'})=H^{p}(Z',({\mathcal{F}}_{2})_{Z'})=0$ for all $p\geq q$, then
by using the long exact sequence of cohomology associated to the exact sequence of sheaves
$$0\longrightarrow Ker (\eta o\psi)\longrightarrow {\mathcal{F}}\stackrel{\eta o\psi}\longrightarrow
({\mathcal{F}}_{2})_{Z'}\longrightarrow 0$$
we see that $H^{p}(X, {\mathcal{F}})=0$ for all $p\geq q$.
\end{proof}
\begin{Thm}{If $\Pi : X\rightarrow \Omega$ is an unbranched Riemann domain and locally
$q$-complete morphism over a Stein space $\Omega$ of dimension $n$,
then $X$ is cohomologically $q$-complete if $n\geq 3$ and $1\leq q\leq n-2$ or if
$n=2$ and $1\leq q\leq 2$}
\end{Thm}
\begin{proof}
Let $\xi : \tilde{\Omega}\rightarrow \Omega$ be a normalization of $\Omega$. If $\tilde{X}$
denotes the fiber product of $\Pi : X\rightarrow \Omega$ and the normalization
$\xi : \tilde{\Omega}\rightarrow \Omega$ , then $\tilde{X}=
\{(x,\tilde{y})\in X\times \tilde{\Omega}: \Pi(x) =\xi(\tilde{y})\}$ and the projection $\Pi_{2} : \tilde{X}\rightarrow \tilde{\Omega}$
is an unbranched Riemann domain over the Stein normal space $\tilde{\Omega}$. ˜
Moreover, since $\Pi_{2}$ is obviously a locally $q$-complete morphism, it follows from lemma $3$ that $\tilde{X}$ is cohomologically $q$-complete. On the other hand, since the projection $\Pi_{1} : \tilde{X}\rightarrow X$ is a finite holomorphic surjection, this implies that X is cohomologically $q$-complete.
\end{proof}
\newpage
\noindent
On behalf of all authors, the corresponding author states that there is no conflict of interest.

\vspace{4pc} \noindent
Youssef Alaoui\\
Hassan II Institute of Agronomy and Veterinary Sciences,\\
Madinat Al Irfane, BP 6202, Rabat, 10101, Morocco,\\
Professor\\
B.P.6202, Rabat-Instituts, 10101. Morocco.\\
Email : {\bf y.alaoui@iav.ac.ma or comp5123ster@gmail.com} \\\\
\end{document}